\documentclass[10pt]{article}
\usepackage{color}
\usepackage{amsmath}
\usepackage{amssymb}
\usepackage{mathrsfs}
\usepackage[utf8]{inputenc}
\usepackage[english]{babel}
\usepackage{amsthm}
\usepackage{tikz}
\usepackage{mathdots}

\usepackage{graphicx}
\graphicspath{{./images/}}

\usepackage[toc,title,page]{appendix}

\usepackage{epsfig}

\usepackage[scaled=1]{helvet}
\usepackage{titlesec}
\usepackage{multido}
\usepackage{multirow}
\usepackage{blindtext}
\usepackage{color}
\usepackage{lipsum}
\usepackage{mathtools}
\usepackage[vcentermath]{youngtab}
\usepackage{young}
\usepackage[all,ps,dvips,graph]{xy}
\usepackage[misc,geometry]{ifsym}

\usepackage[left=3cm,top=2cm,right=3cm,bottom=2cm]{geometry}
\usepackage{graphicx}
\numberwithin{equation}{subsection}

\let\OLDthebibliography\thebibliography
\renewcommand\thebibliography[1]{
  \OLDthebibliography{#1}
  \setlength{\parskip}{0pt}
  \setlength{\itemsep}{0pt plus 0.3ex}
}

\newtheorem{theorem}{Theorem}
\newtheorem{proposition}[theorem]{Proposition}
\newtheorem{corollary}[theorem]{Corollary}
\newtheorem{lemma}[theorem]{Lemma}

\newtheorem{example}[theorem]{Example}
\newtheorem{definition}[theorem]{Definition}

\numberwithin{theorem}{section}

\newcommand{\A}{\mathcal{A}}
\newcommand{\aaa}{\mathfrak{a}}

\newcommand{\Basis}{\mathcal{C}_n}
\newcommand{\BasisOne}{\mathcal{C}^{1,\K}_n}
\newcommand{\B}{\mathcal{B}_{l,m}}
\newcommand{\BB}{\mathcal{B}^{\F}_{l,m}(\kappan,\een)}
\newcommand{\BBnu}{\mathcal{B}^{\F}_{l,m}(\bnu,\een)}
\newcommand{\Bb}{\pmb{\mathfrak{b}}}
\newcommand{\BBB}{\mathcal{B}_{n+1}}
\newcommand{\bbb}{\mathfrak{b}}
\newcommand{\BBk}{\mathcal{B}_k}
\newcommand{\balpha}{\boldsymbol{\alpha}}
\newcommand{\bbeta}{\boldsymbol{\beta}}
\newcommand{\bgamma}{\boldsymbol{\gamma}}
\newcommand{\bGamma}{\boldsymbol{\Gamma}}
\newcommand{\bbkap}{\boldsymbol{\kappa}}
\newcommand\bbn{\mathbf{b}}
\newcommand\bbs{\mathsf{s}}
\newcommand\bbt{\boldsymbol{\mathsf{t}}}
\newcommand\bbu{\mathsf{u}}
\newcommand\bbv{\mathsf{v}}
\newcommand{\Bchico}{\mathcal{B}_{n-1}}
\newcommand\be{\mathbb{E}}
\newcommand\ben{\boldsymbol{e}}
\newcommand{\belongs}{ \in }
\newcommand{\Belongs}{ \ni}
\newcommand{\Bg}{\pmb{\mathfrak{g}}}
\newcommand\bi{\boldsymbol{i}}
\newcommand{\bif}{{\underline{\boldsymbol{i}}}}
\newcommand{\bjf}{{\underline{\boldsymbol{j}}}}
\newcommand\bj{\boldsymbol{j}}
\newcommand\blambda{{\boldsymbol\lambda}}
\newcommand\bmu{{\boldsymbol\mu}}
\newcommand\brho{\boldsymbol{\rho}}
\newcommand\bzeta{{\boldsymbol\zeta}}
\newcommand\bn{\boldsymbol{n}}
\newcommand\bnu{{\boldsymbol\nu}}
\newcommand\boxbluek{\color{blue}\boldsymbol{[k]}}
\newcommand\boxbluej{\color{blue}\boldsymbol{[j]}}
\newcommand\boxredk{\color{red}\boldsymbol{[k]}}
\newcommand\boxredj{\color{blue}\boldsymbol{[j]}}
\newcommand{\bpartial}{\boldsymbol{\partial}}
\newcommand{\bCpartial}{\mathbf{C}\boldsymbol{\partial}}
\newcommand{\Bs}{\pmb{\mathfrak{s}}}
\newcommand\bS{\Sigma}
\newcommand\bs{\mathbf{s}}
\newcommand{\bT}{\pmb{\mathfrak{t}}}
\newcommand\bt{\mathbf{t}}
\newcommand{\bTI}{ \bT^{-1} (\bT_{\theta}^{\blambda}(1))}
\newcommand{\bTII}{ \bT^{-1} (\bT_{\theta}^{\blambda}(2))}
\newcommand{\bTj}{ \bT^{-1} (\bT_{\theta}^{\blambda}(j))}
\newcommand{\bTn}{ \bT^{-1} (\bT_{\theta}^{\blambda}(n))}
\newcommand\btau{{\boldsymbol\tau}}
\newcommand{\Bu}{\pmb{\mathfrak{u}}}
\newcommand{\Bv}{\pmb{\mathfrak{v}}}
\newcommand\bu{\mathbf{u}}
\newcommand\bv{\mathbf{v}}
\newcommand\bcalJ{\boldsymbol{\mathcal{J}}}
\newcommand\bcalL{\boldsymbol{\mathcal{L}}}
\newcommand\bcalm{\boldsymbol{\mathcal{m}}}

\newcommand\bcalNT{\boldsymbol{\mathcal{NT}}}
\newcommand\bcalNGAT{\boldsymbol{\mathcal{NT}}}
\newcommand\bcalD{\boldsymbol{\mathcal{D}}}
\newcommand\bcalR{\boldsymbol{\mathcal{R}}}
\newcommand\bcalV{\boldsymbol{\mathcal{V}}}

\newcommand\bcalU{\boldsymbol{\mathcal{U}}}
\newcommand\bcalY{\boldsymbol{\mathcal{Y}}}
\newcommand\bcalW{\boldsymbol{\mathcal{W}}}
\newcommand\bnabla{\boldsymbol{\nabla}}

\newcommand{\C}{\mathcal{C}_{\rm{YH}}}
\newcommand\calA{\mathcal{A}}
\newcommand\calB{\mathcal{B}}
\newcommand\calC{\mathcal{C}}
\newcommand\calD{\mathcal{D}}
\newcommand\calE{\mathcal{E}}
\newcommand\calF{\mathcal{F}}
\newcommand\calG{\mathcal{G}}
\newcommand\calH{\mathcal{H}}
\newcommand\calL{\mathcal{L}}
\newcommand\calM{\mathcal{M}}
\newcommand\calN{\mathcal{N}}

\newcommand\calO{\mathcal{O}}
\newcommand\calP{\mathcal{P}}
\newcommand\calQ{\mathcal{Q}}
\newcommand\calR{\mathcal{R}}
\newcommand\calS{\mathcal{S}}
\newcommand\calT{\mathcal{T}}

\newcommand\calU{\mathcal{U}}
\newcommand\calV{\mathcal{V}}
\newcommand\calW{\mathcal{W}}
\newcommand\calX{\mathcal{X}}
\newcommand\calY{\mathcal{Y}}
\newcommand\calZ{\mathcal{Z}}

\newcommand{\catorce}{ 14}
\newcommand{\catorceB}{\color{red} 14}
\newcommand{\CC}{ \mathbb C }
\newcommand{\ccc}{\mathfrak{c}}
\newcommand{\ch}{{\rm char}}
\newcommand{\cincuentacinco}{55}
\newcommand{\cincuentacincoR}{\color{red}55}
\newcommand{\Comp}{{\mathcal Comp}_n}
\newcommand{\cuarentacuatro}{ 44}
\newcommand{\cupdot}{\mathbin{\mathaccent\cdot\cup}}
\newcommand{\Cs}{\overleftarrow{C}}
\newcommand{\Cd}{\overrightarrow{C}}

\newcommand{\Der}{{\rm Der}}
\newcommand{\diez}{10}
\newcommand{\dieciseis}{16}
\newcommand{\diezR}{\color{red}10}
\newcommand{\doce}{12}
\newcommand{\doceB}{\color{blue} 12}
\newcommand{\doceR}{\color{red} 12}

\newcommand{\E}{ {\mathcal E}_n(q)}
\newcommand{\e}{\mathfrak{e}}
\newcommand{\EE}{ {\mathcal E}_n}
\newcommand\een{\mathbf{e}}
\newcommand\es{\mathbbm{s}}
\newcommand{\End}{{\rm End}}
\newcommand\et{\mathbbm{t}}

\newcommand\eu{\mathbbm{u}}
\newcommand\ev{\mathbbm{v}}
\newcommand{\Exp}{ {\rm \bf exp} }

\newcommand{\F}{ { \mathbb F}}
\newcommand{\FF}{ {\mathcal F}_n}

\newcommand{\g}{  \mathfrak{g}}
\newcommand{\gl}{\mathfrak{gl}}

\newcommand{\h}{{h}}
\newcommand{\HH}{ \mathcal{H}_n}
\newcommand{\HHO}{ \mathcal{H}^{\OO}_n}
\newcommand{\HHK}{ \mathcal{H}^{\K}_n}
\newcommand{\HHtwo}{ \mathcal{H}_2}
\newcommand{\HHKtwo}{ \mathcal{H}^{\K}_2}
\newcommand{\HHOtwo}{ \mathcal{H}^{\OO}_2}
\newcommand{\HHKOne}{ \mathcal{H}^{1,\K}_n}

\newcommand{\II}{I_{\mathbf{e}}}
\newcommand{\IIa}{I_{e}}
\newcommand{\id}{{\rm id}}
\newcommand{\ind}{{\rm ind}}
\newcommand{\inv}{{\rm inv}}

\newcommand{\JM}{ \mathcal L }

\newcommand{\K}{\mathcal{K}}
\newcommand{\kk}{\mathcal{K}}
\newcommand{\kkk}{k-1}
\newcommand{\kappan}{\boldsymbol{\kappa}}

\newcommand{\Li}{\mathcal{L}}
\newcommand{\LL}{\mathbb{L}}

\newcommand{\m}{\mathfrak{m}}
\newcommand{\MC}{{ {\rm Comp}}_{l,n}}
\newcommand{\MCm}{{ {\rm Comp}}_{l,m}}
\newcommand{\MP}{{\rm Par }_{l,n}}
\newcommand{\mfra}{{\mathfrak{a}}}
\newcommand{\mfrb}{{\mathfrak{b}}}
\newcommand{\mfrc}{{\mathfrak{c}}}
\newcommand{\mfrt}{{\mathfrak{t}}}
\newcommand{\mfrs}{{\mathfrak{s}}}
\newcommand{\mfru}{{\mathfrak{u}}}
\newcommand{\mfrv}{{\mathfrak{v}}}
\newcommand{\mfrw}{{\mathfrak{w}}}
\newcommand{\mfrx}{{\mathfrak{x}}}
\newcommand{\mfry}{{\mathfrak{y}}}
\newcommand{\mfrz}{{\mathfrak{z}}}

\newcommand{\N}{ { \mathbb N}}
\newcommand{\No}{ { \mathbb N}_0}
\newcommand{\nstd}{{\rm NStd}}
\newcommand{\nSLTM}{{$n-$\rm SLTM}}
\newcommand{\nSLTMn}{$\mathcal{TM}_n$}
\newcommand{\NB}{\mathbb{N}\mathcal{B}_k}

\newcommand{\once}{11}
\newcommand{\onceB}{\color{blue}11}
\newcommand{\OnePar}{{ \rm Par}^1_{l,m}}
\newcommand{\OneParn}{{ \rm Par}^1_{l,n}}
\newcommand{\OO}{\mathcal{O}}
\newcommand{\op}{\otimes}

\newcommand{\Par}{{\rm Par}_{l,m}}

\newcommand{\q}{\hat{q}}
\newcommand{\quince}{15}
\newcommand{\quinceB}{\color{red} 15}

\newcommand{\R}{ \mathcal{R}_m}
\newcommand{\Rad}{{\rm Rad}}
\newcommand{\res}{ \textrm{res} }
\newcommand\rrn{\mathbf{r}}

\newcommand{\s}{\mathfrak{s}}
\newcommand{\seq}{{\rm seq}_n}
\newcommand{\shape}{\textsf{shape}}
\newcommand{\Si}{\mathfrak{S}}
\newcommand{\sign}{-}
\newcommand{\sixteenB}{\color{red} 16}
\newcommand{\Snake}{{\rm Snake}}
\newcommand{\spa}{{\rm span}}
\newcommand{\std}{{\rm Std}}
\newcommand{\SLTM}{{\rm SLTM}}
\newcommand{\SLTMn}{$\mathcal{TM}$}

\newcommand{\T}{  \mathfrak{t}}
\newcommand{\tab}{{\rm Tab}}
\newcommand{\tr}{{\rm \textbf{tr}}}
\newcommand{\tR}{ \overline{\R}}
\newcommand{\Tr}{{\rm Tr}}
\newcommand{\trece}{13}
\newcommand{\treintatres}{33}
\newcommand{\treintatresR}{\color{red}33}
\newcommand{\trunc}{ {\B} (\blambda ) }
\newcommand{\truncPrime}{ {\mathbb B}_n^{\prime} (\blambda ) }
\newcommand{\truncSing}{ {\mathbb B}_{\bar n} (\overline{\blambda} ) }
\newcommand{\TT}{{\mathfrak T}}
\newcommand{\TTc}{  \mathcal{T}}

\newcommand{\U}{\mathfrak{u}}
\newcommand{\UU}{\mathbb{U}}
\newcommand{\UUc}{\mathcal{U}}
\newcommand{\ulu}{\underline{u}}
\newcommand{\ulv}{\underline{v}}
\newcommand{\ulw}{\underline{w}}
\newcommand{\ulx}{\underline{x}}
\newcommand{\uly}{\underline{y}}
\newcommand{\ulz}{\underline{z}}

\newcommand{\V}{\mathfrak{v}}
\newcommand{\veintidos}{22}
\newcommand{\VV}{\mathbb{V}}

\newcommand{\Y}{\mathcal{Y}}
\newcommand{\Yy}{\mathcal{Y}_{r,n}}
\newcommand{\yvc}{\Yvcentermath1}
\newcommand{\YY}{\mathcal{Y}_{r,n}(q)}

\newcommand{\Z}{\mathbb{Z}}

\begin{document}
  \Yvcentermath1
\title{On finitude of the number of isomorphism classes in the category  $\bcalNT$}




\author{Diego Lobos Maturana {\thanks{Supported in part by {FONDECYT de Postdoctorado 2024}, N°3240046, and {Subvención a la Instalación en la Academia 2024}, N°85240053, ANID, Chile. }}}

\maketitle


\begin{abstract}
The category $\bcalNT$ was defined in \cite{Lobos2}, it is a category whose objects are commutative nil graded algebras over a field, defined by presentation encoded by triangular matrices. A natural problem related to this category is to reach a complete classification up to isomorphism of its objects. Based in some results coming from \cite{Lobos2}, we can divide this problem by working with encoding matrices of a fixed size $n.$ In \cite{Lobos3} and \cite{Lobos4}, there are several advances for this search, in particular, in \cite{Lobos3} one can see that, for small matrices, the number of isomorphism classes seems to be  finite and independent on the ground field. That fact, opened a series of questions related with the number of isomorphism classes and its relation with the ground field. At that point it was no clear, under which conditions of the ground field, this number could be finite. In this article, among other results, we prove that for each $n\geq4,$ the number of isomorphism classes is finite if and only if the ground field is finite.
\end{abstract}
keywords: Commutative graded algebras. The category $\bcalNT.$
\tableofcontents

\pagenumbering{arabic}

\section{Introduction}
Fix a \emph{ground field} $\F.$ For each strictly lower triangular matrix $T=[t_{ij}]$ of size $n\times n,$ with $n\geq 2$ and entries on $\F,$ we define $\calA(T)$ as the $\F-$algebra given by generators $X_1,\dots,X_n$ and relations:
 \begin{equation}\label{EQ0-def-A(T)}
    \left\{\begin{array}{c}
      X_1^{2}=0 \\
       \quad\\
      X_i^{2}=\sum_{j<i}{t_{ij}}X_jX_i,\quad (2\leq i \leq n)\\
      \quad\\
      X_iX_j-X_jX_i=0,\quad (1\leq i,j \leq n)
    \end{array}\right.
  \end{equation}
 The assignment $\deg(X_i)=2,$ for $1\leq i\leq n,$ converts $\calA(T)$ in a commutative graded algebra.
 The category $\bcalNT=\bcalNT^{{\F}}$ is defined as the category whose objects are all the possible commutative graded algebras $\calA(T)$ defined as above, and whose morphism are all the \emph{preserving degree} homomorphism of graded $\F-$algebras (See Appendix \ref{ssec-prev-def}).

From \cite{Lobos2}, we know that if the matrix $T$ has size $n,$ then the algebra $\calA(T)$ will have dimension equal to $2^{n}.$ This in particular implies that if two encoding matrices $T$ and $S$ have different size, then the algebras $\calA(T)$ and $\calA(S)$ will not be isomorphic. On the other hand, from \cite{Lobos3}, we know that when $T$ and $S$ have the same size, say $n,$ then a necessary and sufficient condition to have an isomorphism $\bgamma:\calA(T)\rightarrow\calA(S),$ is the existence of an invertible $n\times n-$matrix $\Gamma=[\gamma_{ij}]$ whose entries satisfy the following system of equations:

\begin{equation}\label{Eq0-theo-Key-condition}
    2\gamma_{ir}\gamma_{kr}+\gamma_{kr}^{2}s_{ki}=
    \sum_{j<r}{t_{rj}}(\gamma_{kj}\gamma_{kr}s_{ki}
    +\gamma_{kj}\gamma_{ir}+\gamma_{ij}\gamma_{kr})
    ,\quad (1\leq r\leq n;1\leq i<k\leq n).
\end{equation}

In that case, the isomorphism $\bgamma:\calA(T)\rightarrow\calA(S)$ is related with the matrix $\Gamma$ by the equation:
\begin{equation}\label{Eq2-theo-Key-condition}
  \bgamma(X_j)=\sum_{i}\gamma_{ij}Y_i
\end{equation}

where $X_1,\dots,X_n$ are the generators of $\calA(T),$ while $Y_1,\dots,Y_n$ are the generators of $\calA(S)$ (See Appendix \ref{ssec-prev-results}).

Let us denote by $\calT\calM_{n}=\calT\calM_{n}^{\F}$ the set of all strictly lower matrices of size $n\times n$ and entries in $\F.$ We define an equivalence relation on $\calT\calM_{n},$ by stating $T\sim S$ whenever the corresponding algebras $\calA(T)$ and $\calA(S)$ are isomorphic. We denote by $N_n=N_n(\F),$ the number of classes (under relation $\sim$) in $\calT\calM_{n}.$ From \cite{Lobos3}, we know that for $n=2,3$ the number $N_n(\F)$ is finite and independent on the ground field $\F.$  That observation motivated the following questions:

\begin{description}
  \item[Q1:]  Is the number of classes always independent on the ground field?
  \item[Q2:]  Is the number of classes always finite?
\end{description}
 In this article we answer those questions. Propositions \ref{proposition-Nn-not-finite} and \ref{prop2-red-form-n4-W34} imply that in general, the number of classes is not independent on the ground field, moreover under certain conditions this number could be infinite. Since it is clear that the number $N_n(\F)$ is finite whenever the ground field $\F$ is finite, we open the following question:

 \begin{description}
   \item[Q3:] Under which conditions of $\F,$ the number $N_n(\F),$ is finite for all $n$?
\end{description}

The answer to \textbf{Q3} is given in the main result of this article, that is Theorem \ref{theo-Nn-finite}, where we state and prove that the number $N_n(\F)$ is finite if and only if $\F$ is a finite field. Now, knowing the answer to \textbf{Q3}, we open the following question:

 \begin{description}
   \item[Q4:] Let $\F_q$ be a finite field with $q$ elements and let $n$ be an integer greater than 3. Which is the exact value of $N_n(\F_q)$?
\end{description}

For \textbf{Q4}, we only have a partial result. In Theorem \ref{theo-formula-N4}, we provide the explicit value of $N_4(\F_q),$ in function of the cardinality $q$ of the field $\F_q,$ divided in two cases, which depends on if $-1$ is an square in the field $\F_q$ or not. For $n>4,$ the problem turns very complicated, we left it open for future work.

Finally, in \cite{Lobos3}, three \emph{elementary triangular operations} (ETO's for short) were defined:

\begin{equation}\label{eq-ETO-assignments}
  \begin{array}{cc}
    \calP_r(\quad,\alpha):T\mapsto \calP_r(T,\alpha)&(\calP-\textrm{operation} ) \\
    \quad&\quad  \\
     \calF_{(r_1,r_2)}(\quad):T\mapsto \calF_{(r_1,r_2)}(T)& (\calF-\textrm{operation}) \\
    \quad &\quad \\
    \calQ_{(r_0,k_0)}(\quad,\beta):T\mapsto \calQ_{(r_0,k_0)}(T,\beta)& (\calQ-\textrm{operation})
  \end{array}
\end{equation}
each transforming a given matrix $T$ in a new one, whenever $T$ satisfies certain restrictions (see Appendix \ref{ssec-prev-def} for details). Let $T,S$ be two strictly lower triangular $n\times n-$matrices, if there exists of a (finite) sequence of ETO's that transform $T$ into $S,$ then we say that $T$ and $S$ are \emph{equivalent  by ETO's}, and write $T\approx S.$ Equivalence by ETO's is an equivalence relation in $\mathcal{TM}_n.$ Moreover, from \cite{Lobos3} we know that if $T\approx S$, then $T\sim S.$ The converse have not been proved yet but it was conjectured to be true, based on several evidences depicted in  \cite{Lobos3} . That conjecture was called the \textbf{Isomorphism vs ETO's conjecture}.  In Theorem \ref{theo-evidence-ETOs-n4}  we prove that the assertion of the conjecture is true for any pair $T,S$ of strictly lower triangular $4\times 4-$matrices, providing new evidence to support the conjecture.

The content of this paper is organized as follows:
\begin{itemize}
    \item In Section 2, we define \emph{Normalized forms} for the general case, that is a particular set of matrices satisfying certain conditions in their structure. In Proposition \ref{prop-normalized-form}, we prove that any matrix $T\in \mathcal{TM}_n$ is equivalent with a normalized form.
    \item In Subsection \ref{ssec-simplified-form}, we concentrate our attention on $4\times 4-$matrices. For that case we define \emph{Simplified forms}, an specific subset of Normalized forms for $n=4.$ In Proposition  \ref{prop-representatives-n4-general}, we describe a complete set of representatives of classes in $\mathcal{TM}_4$ in terms of simplified forms.
    \item In Propositions \ref{proposition-Nn-not-finite} and \ref{prop2-red-form-n4-W34} we describe explicit criteria to determine when two simplified forms are equivalent or not.
    \item Section 3 is devoted to show some applications of the analysis given in Section 2. In particular in subsection  \ref{ssec-finitude} we obtain the main result of this article, In Theorem \ref{theo-Nn-finite}, we state and prove that the number of classes $N_n(\F)$ is finite if and only if the ground field is finite.
    \item In Subsection \ref{ssec-Explicit-N4} we obtain another important result, in Theorem \ref{theo-formula-N4}, we provide an explicit formula for the number of classes $N_4(\F),$ where  $\F$ is a finite field.
    \item In subsection \ref{ssec-over-C}, we use the algebraic properties and geometrical representation of the field of complex numbers $\mathbb{C}$ to describe explicitly a complete set of representatives of (pairwise different) classes in $\mathcal{TM}_4$ over $\mathbb{C}.$
    \item In subsection \ref{ssec-Iso-vs-ETO-n4}, as an application of the analysis given in Section 2, we provide new evidence to support the \textbf{Isomorphism vs ETO's conjecture} stated in \cite{Lobos3}.
    \item Finally at the end of this article, the reader can find a set of appendices, where we compile notations, definitions and previous results coming from \cite{Lobos2}, \cite{Lobos3} and \cite{Lobos4}.
\end{itemize}

\textbf{Acknowledgments:}

\medskip
It is a pleasure to thank Valentina Delgado and Rahma Salama, for their useful comments on the first draft of this paper.

\medskip
This work was supported by:
\begin{itemize}
  \item {FONDECYT de Postdoctorado 2024}, N°3240046, ANID, Chile.
  \item {Subvención a la Instalación en la Academia 2024}, N°85240053, ANID, Chile.
\end{itemize}


\medskip


\medskip

\section{Normalized and simplified forms}

\subsection{Normalized forms}

In the following we shall denote by $\mathcal{W}_n$ the set of all possible \emph{walls} for matrices $T\in\mathcal{TM}_n.$ Therefore we have:
\begin{equation*}
    \mathcal{W}_{n}=\{(r_1,\dots,r_h): 3\leq r_1\leq \cdots \leq r_h\leq n, h\in[1,n-2]\}\cup\{(0)\}
\end{equation*}
Recall that we define $\mathbb{W}(0_n)=(0).$   Let us denote $\calW_n^{\times}=\calW_n-\{(0)\}.$

For each $\mathbb{W}\in\mathcal{W}_n,$ we define the set $\mathcal{TM}_n(\mathbb{W})$ as the set of all \nSLTM, having $\mathbb{W}$ as a wall. That is,   $\mathcal{TM}_n(\mathbb{W})=\{T\in \mathcal{TM}_n: \mathbb{W}(T)=\mathbb{W}\}.$ We naturally can decompose the set $\mathcal{TM}_n$ into a disjoint union of the sets $\mathcal{TM}_{n}(\mathbb{W})$ with $\mathbb{W}$ running along the set $\calW_n,$ that is:

\begin{equation}\label{eq1-decomp-by-walls}
  \mathcal{TM}_n=\bigcup_{\mathbb{W}\in \mathcal{W}_n}\mathcal{TM}_n(\mathbb{W}).
\end{equation}

If we also define the set $\overline{\mathcal{TM}_n}(\mathbb{W})$ as the set of all the different classes (under relation $\sim$) of matrices in  $\mathcal{TM}_n(\mathbb{W}),$ that is, $\overline{\mathcal{TM}_n}(\mathbb{W})=\{\overline{T}:T\in \mathcal{TM}_n(\mathbb{W})\}.$ By using the decomposition given in equation \ref{eq1-decomp-by-walls}, we can express the number of classes $N_n(\F)$ in $\mathcal{TM}_n$ as the following sum:

\begin{equation}\label{eq2-decomp-by-walls}
N_n(\F)=\sum_{\mathbb{W}\in \mathcal{W}_n}|\overline{\mathcal{TM}_n}(\mathbb{W})|=1+\sum_{\mathbb{W}\in \mathcal{W}_n^{\times}}|\overline{\mathcal{TM}_n}(\mathbb{W})|
\end{equation}

The $1$ at the leftmost expression in equation \ref{eq2-decomp-by-walls} correspond to $|\overline{\mathcal{TM}_n}(0)|,$ since $\overline{\mathcal{TM}_n}(0)=\{\overline{0_n}\}.$

\begin{definition}
    Let $U\nsim 0_n$ be \nSLTM, in $2-REF,$ such that $\mathbb{W}(U)=(r_1,\dots,r_h).$ We say that $U$ is \emph{normalized}, if the following conditions hold
    \begin{itemize}
        \item If $u_{r_1,j}\neq 0,$ for some  $j\in[1,r_1[,$ then $u_{r_1,j}=1.$
        \item For each $r\in [r_1,n],$ the leader $u_{rc}$ of the $r-$th row, is equal to $1.$
    \end{itemize}

    We also say that the zero matrix is \emph{normalized}
\end{definition}

\begin{proposition}\label{prop-normalized-form}
    Let $T$ be a \nSLTM, then $T$ is equivalent by ETO's to a normalized matrix.
\end{proposition}

\begin{proof}
    We already know from \cite{Lobos4}, that $T$ is equivalent by ETO's to a matrix $S$ that is in the $2-$REF. If $T\nsim 0_n,$ we shall prove that for each matrix $S$ in $2-$REF, we can obtain a matrix $U$ via a finite sequence of $\calP-$operations applied (progressively) to $S:$
    \begin{enumerate}
        \item  For each $j\in [1,r_1[,$ such that $s_{rj}\neq 0,$ we apply the operation $\calP_{j}(\quad,1/s_{rj}).$ The obtained matrix after this process, say $V,$ clearly satisfy the condition that for each $j\in[1,r_1[,$ such that $v_{r_1j}\neq 0,$ we necessarily have that $v_{r_1j}=1.$
        \item By our construction we have that the leader $v_{r_1c_1}$ of the $r-$th row is equal to $1.$ Assuming that $r\in ]r_1,n]$ and for each $r'\in [r_1,r[$  the leader $v_{r'c'}$ of the $r'-$th row is equal to $1.$ If $v_{rc}$ is the leader of the $r-$th row of $V,$ then we apply the operation $\calP_{r}(\quad,v_{rc})$ obtaining a new matrix, say $W,$ that is in the $2-$REF, such that $W^{[r_1]}=V^{[r_1]}$ and for which  for each $r'\in [r_1,r]$ the leader $w_{r'c'}$ of the $r'-$th row, is equal to $1.$
        \item By a recursive argument, we obtain a new matrix $U$ in the $2-$REF, normalized and equivalent by ETO's to $T.$
    \end{enumerate}
\end{proof}
\begin{example}\label{example-normalized-form}
Assuming $\F$ is a field of characteristic zero, consider the matrix
\begin{equation*}
    T=\left[\begin{matrix}
        0&0&0&0\\
        0&0&0&0\\
        2&1&0&0\\
        12&9&3&0
    \end{matrix}\right]
\end{equation*}
One can easily check that $T$ is a matrix in $2-$REF, such that $\mathbb{W}(T)=(3,4).$ This matrix is not normalized since, for example, the leader of the $4-$th row is not equal to $1.$ If we follow the algorithm described in the proof of proposition \ref{prop-normalized-form}, we can see that:
\begin{equation}\label{eq-example-normalized-form}
    T\approx V:=\calP_1(T,1/2)= \left[\begin{matrix}
        0&0&0&0\\
        0&0&0&0\\
        1&1&0&0\\
        6&9&3&0
    \end{matrix}\right]\approx U:= \calP_4(V,3)= \left[\begin{matrix}
        0&0&0&0\\
        0&0&0&0\\
        1&1&0&0\\
        2&3&1&0
    \end{matrix}\right]
\end{equation}
The matrix $U$ of equation \ref{eq-example-normalized-form} is normalized, as desired.
\end{example}
As a consequence of Proposition \ref{prop-normalized-form}, we have the following:

\begin{corollary}\label{coro-normalized-forms}
Let $\mathbb{W}\in \calW_n^{\times}.$ Then we have
    \begin{equation*}
        \overline{\mathcal{TM}_n}(\mathbb{W})=\left\{\overline{U}:U\in \mathcal{TM}_n(\mathbb{W})\quad \textrm{and}\quad U\quad \textrm{is normalized}\right\}.
    \end{equation*}
\end{corollary}
Corollary \ref{coro-normalized-forms} basically says that if we want to obtain a complete set of representatives of classes for $\overline{\mathcal{TM}_n}(\mathbb{W}),$ we can focus our attention in normalized forms. Corollary \ref{coro-normalized-forms} does not imply that two normalized forms are not equivalents, so if we want to obtain a complete set of representatives of pairwise classes, we should look for necessary and sufficient conditions to determine when two normalized forms are equivalent or not.  In the following subsection we do so for the case $n=4.$

\subsection{Simplified forms in $\mathcal{TM}_4$}\label{ssec-simplified-form}

For the rest of this section we consider the following notation:

\begin{equation}\label{eq-notation-mat-B}
     0_4=\left[\begin{matrix}
        0&0&0&0\\
        0&0&0&0\\
        0&0&0&0\\
        0&0&0&0
    \end{matrix}\right];\quad
    B=\left[\begin{matrix}
        0&0&0&0\\
        0&0&0&0\\
        0&0&0&0\\
        1&1&0&0
    \end{matrix}\right];\quad B'=\left[\begin{matrix}
        0&0&0&0\\
        0&0&0&0\\
        0&0&0&0\\
        1&1&1&0
    \end{matrix}\right]
\end{equation}

\begin{equation}\label{eq-notation-mat-C-D}
 C_{a}=\left[\begin{matrix}
        0&0&0&0\\
        0&0&0&0\\
        1&1&0&0\\
        a&0&1&0
    \end{matrix}\right],\quad(a\in\F);\quad  D_{u}=\left[\begin{matrix}
            0&0&0&0\\
            0&0&0&0\\
            1&1&0&0\\
            u&1&0&0
        \end{matrix}\right],\quad(u\in \F^{\times}).
\end{equation}

\begin{proposition}\label{prop-representatives-n4-general}
   Independently on the ground field, we have:
    \begin{equation}\label{eq1-prop-representatives-n4-general}
    \begin{array}{cc}
         \overline{\mathcal{TM}_4}(0)=\{\overline{0}_4\},& \overline{\mathcal{TM}_4}(4)=\{\overline{B},\overline{B'}\},  \\
         \quad& \quad \\
         \overline{\mathcal{TM}_4}(3,4)=\{\overline{C}_a:a\in \F\},& \overline{\mathcal{TM}_4}(3)=\{\overline{D}_u:u\in \F^{\times}\}
    \end{array}
    \end{equation}
    \end{proposition}
\begin{proof}
As we already mentioned the matrix $0_4$ is the only normalized matrix in the corresponding sets $\calT\calM_{4}(0).$ On the other hand, following the proof of Proposition \ref{prop-normalized-form}, we can see that for any matrix $T\in \calT\calM_4(4),$ we have that:
\begin{itemize}
    \item
\begin{equation}
    T\sim B,\quad \textrm{if}\quad \mathbb{M}(T)=\left(\left[\begin{matrix}
        4\\2
    \end{matrix}\right]\right)
\end{equation}
\item
\begin{equation}
    T\sim B',\quad \textrm{if}\quad \mathbb{M}(T)=\left(\left[\begin{matrix}
        4\\3
    \end{matrix}\right]\right)
\end{equation}
\end{itemize}
(See appendix \ref{ssec-prev-def} for definition of $\mathbb{M}(T)$).

For the case of a matrix $T\in \calT\calM_4(3),$ if we follow the proof of Proposition \ref{prop-normalized-form},  we can see that $T\sim D_u$ for some $u\in \F^{\times}.$

Finally if $T\in \calT\calM_4(3,4),$ if we follow the proof of Proposition \ref{prop-normalized-form},  we can see that $T$ is equivalent (by ETO's) to a normalized matrix of the form:
\begin{equation}\label{eq-U-normal}
    V=\left[\begin{matrix}
        0&0&0&0\\
        0&0&0&0\\
        1&1&0&0\\
        u&v&1&0
    \end{matrix}\right]
\end{equation}
The proof that $V$ can be reduced to a matrix of the form $C_a$ is given in Lemma \ref{lemma1-red-form-n4-W34} below.
\end{proof}

\begin{lemma}\label{lemma1-red-form-n4-W34}
Given $u,v\in\F,$ consider the normalized matrix $V$ of equation \ref{eq-U-normal}.
If $a=2uv+u+v,$ then the matrices $V$ and $C_a$ are equivalent by ETO's.
\end{lemma}
\begin{proof}
Take $\epsilon\in \F^{\times}$ and define $T_1=\calQ_{(2,1)}(V,\epsilon).$ Then $T_1$ is well-defined and looks like this:
\begin{equation}
    T_1=\left[\begin{matrix}
        0&0&0&0\\
        -2\epsilon&0&0&0\\
        1+\epsilon&1&0&0\\
        u+v\epsilon&v&1&0
    \end{matrix}\right],
\end{equation}
let $\beta:=\frac{\epsilon-1}{2\epsilon}.$ Note that $\Delta_{(1,2,3)}^{(2)}(T_1)=-2\epsilon \beta,$ therefore the matrix $T_2:=\calQ_{(3,2)}(T_1,\beta)$ is well-defined and looks like this:
\begin{equation}
    T_2=\left[\begin{matrix}
        0&0&0&0\\
        -2\epsilon&0&0&0\\
        1+\epsilon&1/\epsilon&0&0\\
        u+v\epsilon&v'/(2\epsilon)&1&0
    \end{matrix}\right],\quad\textrm{where}\quad v'={(2v+1)\epsilon-1}.
\end{equation}
Now let $T_3:=\calP_{2}(T_2,\epsilon),$ then
\begin{equation}
    T_3=\left[\begin{matrix}
        0&0&0&0\\
        -2&0&0&0\\
        1+\epsilon&1&0&0\\
        u+v\epsilon&v'/2&1&0
    \end{matrix}\right].
\end{equation}
Let $T_4:=\calQ_{(2,1)}(T_3,-1),$ then it is clear that $T_4$ is well-defined and looks like this:
\begin{equation}
    T_4=\left[\begin{matrix}
        0&0&0&0\\
        0&0&0&0\\
        \epsilon&1&0&0\\
        u' &v'/2&1&0
    \end{matrix}\right],\quad \textrm{where}\quad u'=\frac{2u+1-\epsilon}{2}
\end{equation}
Now let $T_5=\calP_{1}(1/\epsilon),$ then we have
\begin{equation}
    T_5=\left[\begin{matrix}
        0&0&0&0\\
        0&0&0&0\\
        1&1&0&0\\
        u'/\epsilon &v'/2&1&0
    \end{matrix}\right].
\end{equation}
Finally, define $U_{\epsilon}=\calF_{(1,2)}(T_5).$ It is easy to see that $U_{\epsilon}$ is well defined and looks like this:
\begin{equation}
    U_{\epsilon}=\left[\begin{matrix}
        0&0&0&0\\
        0&0&0&0\\
        1&1&0&0\\
        v'/2& u'/\epsilon &1&0
    \end{matrix}\right].
\end{equation}
By construction we have $V$ is equivalent by ETO's to $U_{\epsilon}.$ It is not difficult to see that when we replace $\epsilon:=2u+1,$ we obtain $U_{\epsilon}=C_a.$ Therefore $V$ is equivalent by ETO's to $C_a$ as desired.
\end{proof}
A matrix of the form $0_4,B,B',D_u$ or $C_a$ (as in equations \ref{eq-notation-mat-B} and \ref{eq-notation-mat-C-D}) will be said that is in a \emph{simplified form}. Proposition \ref{prop-normalized-form} and Lemma \ref{lemma1-red-form-n4-W34}, basically imply that each normalized matrix in $\mathcal{MT}_4$, can be reduced (by an ETO's process) into an adequate simplified form. Moreover Proposition \ref{prop-representatives-n4-general} imply that we can choose a complete set of representatives of classes for $\mathcal{MT}_4$, composed only by simplified forms.
\begin{example}
A simplified form equivalent by ETO's to the normalized matrix $U$ of equation \ref{eq-example-normalized-form} is $C_{17}.$ That is
    \begin{equation*}
        \left[\begin{matrix}
        0&0&0&0\\
        0&0&0&0\\
        1&1&0&0\\
        2&3&1&0
    \end{matrix}\right]\approx  \left[\begin{matrix}
        0&0&0&0\\
        0&0&0&0\\
        1&1&0&0\\
        17&0&1&0
    \end{matrix}\right]
    \end{equation*}
\end{example}

 As a direct consequence of Proposition \ref{prop-representatives-n4-general}, we have:

\begin{corollary}\label{coro1-rep-n4}
    Independently on the ground field, the number $N_4(\F)$ is given by
    \begin{equation}\label{eq-coro-rep-n4}
        N_4(\F)=3+|\{\overline{D_u}:u\in \F^{\times}\}|+|\{\overline{C_a}:a\in \F|\}.
    \end{equation}
    In particular if $\F=\F_q$ is a finite field with $q$ elements (and $\textrm{char}(\F)\neq 2$), then $N_4(\F_q)\leq 2q+2.$
\end{corollary}


Equation  $\overline{\mathcal{TM}_4}(3,4)=\{\overline{C}_a:a\in \F\},$ (resp.  $ \overline{\mathcal{TM}_4}(3)=\{\overline{D}_a:a\in \F^{\times}\}$) in Proposition \ref{prop-representatives-n4-general}, does not imply that the set $\{\overline{C}_a:a\in \F\}$ (resp. $\{\overline{D}_a:a\in \F^{\times}\}$) is a complete set of pairwise different classes for $\mathcal{TM}_4(3,4)$ (resp. for $\mathcal{TM}_4(3)$). Therefore, in order to obtain complete sets of representative of pairwise different classes for $\mathcal{TM}_4(3,4)$ and for $\mathcal{TM}_4(3)$, we need to find more specific criteria. In propositions \ref{proposition-Nn-not-finite} and \ref{prop2-red-form-n4-W34} below, we obtain necessary and sufficient conditions to determinate if two given simplified forms are equivalent or not.


\begin{proposition}\label{proposition-Nn-not-finite}
  Let $u,v\in \F^{\times}$ be such that  $u\neq v.$ Then $D_u\sim D_v$ if and only if $\frac{u}{v}$ is an square in $\F.$
\end{proposition}

\begin{proof}
Lets assume that $\bgamma:\calA(D_u)\rightarrow\calA(D_v)$ is an isomorphisms with associated matrix $\Gamma=[\gamma_{ij}].$ By Proposition \ref{prop-Gamma-triang-2}, we know that $\Gamma$ have one of the following forms:
   \begin{equation}\label{eq1-proof-proposition-Nn-not-finite}
     \Gamma_1=\left[\begin{matrix}
                    a & 0 & x & z \\
                    0 & b & y & w \\
                    0 & 0 & c & 0 \\
                    0 & 0 & 0 & d
                  \end{matrix}\right],
     \Gamma_2=\left[\begin{matrix}
                    0 & a & x & z \\
                    b & 0 & y & w \\
                    0 & 0 & c & 0 \\
                    0 & 0 & 0 & d
                  \end{matrix}\right],
     \Gamma_3=\left[\begin{matrix}
                    a & 0 & z & x \\
                    0 & b & w & y \\
                    0 & 0 & 0 & c \\
                    0 & 0 & d & 0
                  \end{matrix}\right],
   \Gamma_4=\left[\begin{matrix}
                    0 & a & z & x \\
                    b & 0 & w & y \\
                    0 & 0 & 0 & c \\
                    0 & 0 & d & 0
                  \end{matrix}\right].
   \end{equation}
   where $abcd\neq0.$ Also, we know that the entries of $\Gamma$ are related by equation \ref{Eq1-theo-Key-condition} of Theorem \ref{theo-Key-condition}.

   We shall assume that $\Gamma$ has the form $\Gamma_1$ in equation \ref{eq1-proof-proposition-Nn-not-finite}, the analysis are equivalent for all the forms. Therefore if we replace $r=3$ in equation \ref{Eq1-theo-Key-condition}, we obtain (after simplifications), the following system of equations:

   \begin{equation}\label{eq2-proof-proposition-Nn-not-finite}
     \left\{\begin{array}{c}
              2xy=ay+bx \\
              2x=a-c \\
              2y=b-c
            \end{array}\right.
   \end{equation}
   Combining those equations, one obtain
   \begin{equation}\label{eq3-proof-proposition-Nn-not-finite}
   c^2=ab
   \end{equation}
  Therefore, the parameters $a,b$ have to be taken such that $ab$ is an square in the field $\F.$ If we do so, any choice of $c=\pm \sqrt{ab}$ define $x$ and $y$ and we can check that $\gamma(X_3)=xY_1+yY_2+cY_3$ satisfies the desired conditions.

  Now if we replace $r=4$ in equation \ref{Eq1-theo-Key-condition}, we obtain (after simplifications), the following system of equations:
  \begin{equation}\label{eq4-proof-proposition-Nn-not-finite}
     \left\{\begin{array}{c}
              2zw=uaw+bz \\
              2z=ua-vd \\
              2w=b-d
            \end{array}\right.
  \end{equation}
  Combining those equations we obtain
  \begin{equation}\label{eq5-proof-proposition-Nn-not-finite}
   vd^2=uab
   \end{equation}
   from where we conclude that, for having $D_u\sim D_v$ it is necessary that the fraction $\frac{u}{v}$ is an square on $\F.$

   On the other hand, if we assume that $\frac{u}{v}$ is an square on $\F,$ say $\frac{u}{v}=\epsilon^2,$ for some $\epsilon\in \F^{\times},$ we can check that the matrix
   \begin{equation}\label{eq6-proof-proposition-Nn-not-finite}
     \Gamma=\left[\begin{matrix}
                    1 & 0 & 0 & \frac{1}{2}(u-v\epsilon) \\
                    \quad&\quad&\quad&\quad\\
                    0 & 1 & 0 & \frac{1}{2}(1-\epsilon) \\
                    \quad&\quad&\quad&\quad\\
                    0 & 0 & 1 & 0 \\
                    \quad&\quad&\quad&\quad\\
                    0 & 0 & 0 & \epsilon
                  \end{matrix}\right]
   \end{equation}
   defines an isomorphism $\bgamma:\calA(D_u)\rightarrow \calA(D_v).$
\end{proof}


\begin{proposition}\label{prop2-red-form-n4-W34}
    Let $a,a'\in \F$ with $a\neq a'.$ Then $C_{a}\sim C_{a'}$ if and only if $2a+1\in(\F^{\times})^{2}$ and $a'=-\frac{a}{2a+1}.$
\end{proposition}
\begin{proof}
    Suppose that $2a+1\in(\F^{\times})^{2}$ and  $a'=-\frac{a}{2a+1}.$ Let $u=a$ and $v=0,$ then by the proof of Lemma \ref{lemma1-red-form-n4-W34}, we have that
    \begin{equation}
        T=\left[\begin{matrix}
        0&0&0&0\\
        0&0&0&0\\
        1&1&0&0\\
        u&v&1&0
    \end{matrix}\right]\approx T'=\left[\begin{matrix}
        0&0&0&0\\
        0&0&0&0\\
        1&1&0&0\\
        u'/\epsilon&v'/2&1&0
    \end{matrix}\right]
    \end{equation}
    where $u'=\frac{2a+1-\epsilon}{2}$ and $v'=\epsilon-1,$ for some $\epsilon\in \F^{\times}$. If we take $\epsilon=\sqrt{2a+1},$ we obtain
\begin{equation}
        T=\left[\begin{matrix}
        0&0&0&0\\
        0&0&0&0\\
        1&1&0&0\\
        u&v&1&0
    \end{matrix}\right]\approx \left[\begin{matrix}
        0&0&0&0\\
        0&0&0&0\\
        1&1&0&0\\
        u''& v''&1&0
    \end{matrix}\right]
    \end{equation}
where $u''=v''=\frac{\sqrt{2a+1}-1}{2}.$
By Lemma \ref{lemma1-red-form-n4-W34}, we have that
\begin{equation}
      \left[\begin{matrix}
        0&0&0&0\\
        0&0&0&0\\
        1&1&0&0\\
        u''& v''&1&0
    \end{matrix}\right]\approx
    \left[\begin{matrix}
        0&0&0&0\\
        0&0&0&0\\
        1&1&0&0\\
        a''& 0&1&0
    \end{matrix}\right]
    \end{equation}
    where $a''=2u''v''+u''+v''.$ After some calculus, we can check that $a''=a'=-\frac{a}{2a+1}.$ and therefore $C_{a}\approx C_{a'}.$ From \cite{Lobos3}, we know that this imply that $C_{a}\sim C_{a'}$ (see Theorem \ref{theo-sim-vs-approx}).

    On the other hand, if we assume that $C_{a}\sim C_{a'}.$ then there is an invertible matrix $\Gamma=[\gamma_{ij}]$ whose entries satisfy equation \ref{Eq1-theo-Key-condition} (Taking $T=C_{a}, S=C_{a'}$). By Proposition \ref{prop-Gamma-triang-2}, we also know that the matrix $\Gamma$ has one of the following forms
    \begin{equation}
        \Gamma_1=\left[\begin{matrix}
            {\color{blue}\gamma_{11}}&{\color{blue}0}&\gamma_{13}&\gamma_{14}\\
            {\color{blue}0}&{\color{blue}\gamma_{22}}&\gamma_{23}&\gamma_{24}\\
            0&0&{\color{red}\gamma_{33}}&\gamma_{34}\\
            0&0&0&\boldsymbol{\gamma_{44}}
        \end{matrix}\right]\quad\textrm{or}\quad
        \Gamma_2=\left[\begin{matrix}
            {\color{blue}0}&{\color{blue}\gamma_{12}}&\gamma_{13}&\gamma_{14}\\
            {\color{blue}\gamma_{21}}&{\color{blue}0}&\gamma_{23}&\gamma_{24}\\
            0&0&{\color{red}\gamma_{33}}&\gamma_{34}\\
            0&0&0&\boldsymbol{\gamma_{44}}
        \end{matrix}\right]
    \end{equation}
    we can assume with no losing of generality, that $\Gamma$ has the form of $\Gamma_1$ (the analysis is equivalent for $\Gamma_2$). Since $\det(\Gamma)\neq 0,$ we have that $\gamma_{ii}\neq 0$ for each $i\in[1,4].$ Using equation \ref{Eq1-theo-Key-condition}, we deduce the following system of equations:
    \begin{equation}\label{system-aa}
        \left\{\begin{array}{cc}
               2\gamma_{13}\gamma_{23}=&\gamma_{11}\gamma_{23}+\gamma_{22}\gamma_{13}\\

             2\gamma_{13}=&\gamma_{11}-\gamma_{33},  \\

              2\gamma_{23}=&\gamma_{22}-\gamma_{33},  \\

              2\gamma_{14}\gamma_{24}=&a\gamma_{11}\gamma_{24}+\gamma_{23}\gamma_{14}+\gamma_{13}\gamma_{24}\\

             2\gamma_{14}=&a\gamma_{11}-a'\gamma_{44}+\gamma_{13}\\

             2\gamma_{24}=&\gamma_{23}\\

             2\gamma_{34}=&\gamma_{33}-\gamma_{44}\\

             \gamma_{34}^{2}=&\left(a\gamma_{11}+\gamma_{33}+\gamma_{13}-2\gamma_{14} \right)\gamma_{34}+\gamma_{33}\gamma_{14}\\
             \gamma_{34}^{2}=&\left(\gamma_{33}+\gamma_{23}-2\gamma_{24}\right)\gamma_{34}+\gamma_{33}\gamma_{24}
        \end{array}\right.
    \end{equation}

   Combining equations in \ref{system-aa}, we obtain:

     \begin{equation}\label{system-ab}
        \left\{\begin{array}{cc}
              \gamma_{33}^{2}=&\gamma_{11}\gamma_{22}\\

             \gamma_{44}^{2}=&\gamma_{22}\gamma_{33}\\

             (2a'+1)\gamma_{44}^{2}=&(2a+1)\gamma_{11}\gamma_{33}\\
              0=&\left(\gamma_{22}-\gamma_{33}\right)\left((2a+1)\gamma_{11}-\gamma_{33}\right)
        \end{array}\right.
    \end{equation}

    In particular we can see that $2a+1=0$ if and only if $2a'+1=0,$ but this only happens when $a=a'=-1/2.$ By our assumption $a\neq a',$ then we conclude that $2a+1\neq 0$ and $2a'+1\neq 0.$ Using that fact, we can see that
    \begin{equation}\label{eq1-post-system-aa}
        \frac{\gamma_{11}}{\gamma_{22}}=\frac{2a'+1}{2a+1}
    \end{equation}

    On the other hand, the last equation in \ref{system-ab}, we have that $\gamma_{22}=\gamma_{33}$ or $(2a+1)\gamma_{11}=\gamma_{33}.$ But if we assume that $\gamma_{22}=\gamma_{33},$ we also obtain that $\gamma_{11}=\gamma_{22}=\gamma_{33},$ then equation \ref{eq1-post-system-aa} implies that $a=a',$ and this contradicts our assumption. Therefore we conclude that $(2a+1)\gamma_{11}=\gamma_{33}.$

    Since $\gamma_{44}^{2}=\gamma_{22}\gamma_{33}=(2a+1)\gamma_{11}\gamma_{22}=(2a+1)\gamma_{33}^{2},$ we conclude that a necessary condition for having $T\sim S$ is that $2a+1\in(\F^{\times})^{2}.$ On the other hand, since $(2a'+1)\gamma_{44}^{2}=(2a+1)\gamma_{11}\gamma_{33}=\gamma_{33}^{2},$ one conclude that $2a'+1=1/(2a+1),$ that is, $a'=-a/(2a+1),$ as desired.
\end{proof}

\begin{corollary}\label{coro-prop2-red-form-n4-W34}
    Suppose that $a\in \F$ is such that for any $a'\in \F,$ the condition $C_{a}\sim C_{a'}$ implies that $a=a'.$ Then $a\in\{x\in\F: 2x+1\notin (\F^{\times})^{2}\}\cup \{0,-1,-1/2\}.$
\end{corollary}

\begin{proof}
    It follows directly from Proposition \ref{prop2-red-form-n4-W34} and the fact that the only two solutions of the equation $a=-a/(2a+1)$ are $a=0$ and $a=-1.$ For the particular case of $a=-1/2,$ it is easy to see that $2a+1=0\notin (\F^{\times})^{2}.$
\end{proof}

\section{Applications}


\subsection{On finitude of $N_n(\F)$}\label{ssec-finitude}


Let $a\in\F$ and $n\in\mathbb{N}.$ For the remainder of this section, we shall denote by $C_{a,n}=[c_{ij}]$ the \nSLTM, whose entries are given by:

\begin{equation}
    c_{ij}=\left\{\begin{array}{cc}
         1& \textrm{if}\quad i\in[3,n-1],j\in[1,2] \\
         1& \textrm{if}\quad (i,j)=(n,3)\\
         a& \textrm{if}\quad (i,j)=(n,1)\\
         0& \textrm{otherwise}
    \end{array}\right.
\end{equation}

For example
\begin{equation*}
 C_{a,4}=C_a=\left[\begin{matrix}
        0&0&0&0\\
        0&0&0&0\\
        1&1&0&0\\
        a&0&1&0
    \end{matrix}\right],\quad C_{a,5}=\left[\begin{matrix}
        0&0&0&0&0\\
        0&0&0&0&0\\
        1&1&0&0&0\\
        1&1&0&0&0\\
        a&0&1&0&0
    \end{matrix}\right]
\end{equation*}

\begin{proposition}\label{prop3-red-form-n-W3n}
    Let $n>4,$ $a,a'\in \F$ with $a\neq a'.$ A necessary condition for having $C_{a,n}\sim C_{a',n}$ is that $2a+1\in\left(\F^{\times}\right)^{2}$ and $a'=-\frac{a}{2a+1}.$
\end{proposition}

\begin{proof}
    Assuming $C_{a,n}\sim C_{a',n},$ then there is an invertible matrix $\Gamma=[\gamma_{ij}]$ satisfying equation \ref{Eq1-theo-Key-condition}, taking $T=C_{a,n}, S=C_{a',n}.$ By proposition \ref{prop-Gamma-triang-2}, we can assume with no loss of generality that $\Gamma$ has the following form:

    \begin{equation}
    \Gamma=\left[\begin{matrix}
        {\color{blue}\gamma_{11}}& {\color{blue}0}&\gamma_{13}&\gamma_{14}&\cdots &\gamma_{1,n-1}&\gamma_{1n}\\
         {\color{blue}0}&  {\color{blue}\gamma_{21}}& \gamma_{23}&\gamma_{24}&\cdots& \gamma_{1,n-1}&\gamma_{2n}\\
        0&0 &        {\color{red}\gamma_{33}}& {\color{red}0}&\cdots &{\color{red}0}&\gamma_{3n}\\
        0&0&{\color{red}0}& {\color{red}\gamma_{44}}&\cdots&{\color{red}0}&\gamma_{4n}\\
        \vdots& \vdots& \vdots&\quad&\ddots &\vdots&\vdots\\
        0&0&{\color{red}0}&{\color{red}0}&\cdots&{\color{red}\gamma_{n-1,n-1}}&\gamma_{n-1,n}\\
        0&0&0&0&\cdots&0&\boldsymbol{\gamma_{nn}}
    \end{matrix}\right]
    \end{equation}

     Since $\det(\Gamma)\neq 0,$ we have that $\gamma_{ii}\neq 0$ for each $i\in[1,n].$ Using equation \ref{Eq1-theo-Key-condition}, we deduce the following system of equations:
    \begin{equation}\label{system-ac}
        \left\{\begin{array}{ccc}
               2\gamma_{1r}\gamma_{2r}=&\gamma_{11}\gamma_{2r}+\gamma_{2r}\gamma_{1r},&\textrm{for}\quad r\in[3,n-1]\\

             2\gamma_{1r}=&\gamma_{11}-\gamma_{rr},&\textrm{for} \quad r\in[3,n-1] \\

              2\gamma_{2r}=&\gamma_{22}-\gamma_{rr},&\textrm{for} \quad r\in[3,n-1] \\

              2\gamma_{1n}\gamma_{2n}=&a\gamma_{11}\gamma_{2n}+\gamma_{23}\gamma_{1n}+\gamma_{13}\gamma_{2n},&\quad\\

             2\gamma_{1n}=&a\gamma_{11}-a'\gamma_{nn}+\gamma_{13},&\quad \\

             2\gamma_{2n}=&\gamma_{23},&\quad \\

             2\gamma_{3n}=&\gamma_{33}-\gamma_{nn},&\quad \\

             \gamma_{3n}^{2}=&\left(a\gamma_{11}+\gamma_{33}+\gamma_{13}-2\gamma_{1n} \right)\gamma_{3n}+\gamma_{33}\gamma_{1n},&\quad\\
             \gamma_{3n}^{2}=&\left(\gamma_{33}+\gamma_{23}-2\gamma_{2n}\right)\gamma_{3n}+\gamma_{33}\gamma_{2n},&\quad
        \end{array}\right.
    \end{equation}
    (compare with equation \ref{system-aa}). An analogous analysis as we did in the proof of Proposition \ref{prop2-red-form-n4-W34}, implies that a necessary condition for having $C_{a,n}\sim C_{a',n}$ is that $2a+1\in\left(\F^{\times}\right)^{2}$ and $a'=-\frac{a}{2a+1}.$
\end{proof}

Note that the system \ref{system-ac} used in the proof of Proposition \ref{prop3-red-form-n-W3n}, was obtained taking some of the equations coming from equation \ref{Eq1-theo-Key-condition}, but not all of them, therefore we cannot claim that the converse is true for $n>4.$

\begin{theorem}\label{theo-Nn-finite}
The number $N_n(\F)$ is finite if and only if $\F$ is a finite field.
\end{theorem}

\begin{proof}
It is clear that if $\F$ is finite, then $N_n(\F)$ is finite. We shall prove that $N_n(\F)$ is infinite whenever $\F$ is an infinite field:

Let us choose an arbitrary element $a_1\in \F,$ we define the set:
\begin{equation*}
    A_1=\left\{\begin{array}{cc}
         \left\{a_1,\frac{-a_1}{2a_1+1}\right\},&\textrm{if}\quad 2a_1+1\in \left(\F^{\times}\right)^{2} \\
         \quad& \quad\\
         \{a_1\},&\textrm{otherwise}
    \end{array}\right.
\end{equation*}
Since $\F$ is infinite and $A_1\subset \F$ is a finite set, we can choose $a_2\in \F-A_1$ and define the set:
\begin{equation*}
    A_2=\left\{\begin{array}{cc}
         A_1\cup \left\{a_2,\frac{-a_2}{2a_2+1}\right\},&\textrm{if}\quad 2a_2+1\in \left(\F^{\times}\right)^{2}  \\
         \quad& \quad\\
         A_1\cup\{a_2\},&\textrm{otherwise}
    \end{array}\right.
\end{equation*}
Following recursively, for $k>2$ assume that we have already defined the finite subset $A_{k-1}\subset\F$ and choose an element $a_{k}\in \F-A_{k-1}$ (this is possible, since $\F$ is an infinite field). Now define the set:
\begin{equation*}
    A_k=\left\{\begin{array}{cc}
         A_{k-1}\cup \left\{a_k,\frac{-a_k}{2a_k+1}\right\},&\textrm{if}\quad 2a_k+1\in \left(\F^{\times}\right)^{2}  \\
         \quad& \quad\\
         A_{k-1}\cup\{a_k\},&\textrm{otherwise}
    \end{array}\right.
\end{equation*}
Therefore, we have chosen a sequence $(a_k)_{k\in \mathbb{N}}$ of not repeated elements of $\F.$ By Propositions \ref{prop2-red-form-n4-W34} and \ref{prop3-red-form-n-W3n} we conclude that our construction implies that the set of matrices $ \left\{C_{a_k,n}:k\in \mathbb{N}\right\},$ is an infinite set of representative of pairwise different classes (under relation $\sim$).
\end{proof}

\subsection{Explicit formula for $N_4(\F_q).$ }\label{ssec-Explicit-N4}
In this subsection we consider as a ground field $\F=\F_{q},$ the finite field of $q$ elements (with $\textrm{char}(\F)\neq 2$). We shall assume certain familiarity with the theory of finite fields. For some previous background on finite fields the reader can see \cite{Lidl-finitefields} for example.



\begin{lemma}\label{lemma2-count-N4-q}
    If $\F$ is a finite field with characteristic different from $2,$ then $|\overline{\mathcal{TM}_4}(3)|=2.$
\end{lemma}

\begin{proof}
Let $\F=\F_q$ be a finite set with $q$ elements and characteristic different from $2.$
 Let $T\in \calM_4(3),$ then by Proposition \ref{prop-representatives-n4-general}, we have that   $T\sim D_a,$ for some $a\in\F^{\times}.$ We also know from Proposition \ref{proposition-Nn-not-finite}, that for any pair $a,b\in\F^{\times}$ we have that $D_a\sim D_b$ if and only if $b/a\in (\F^{\times})^2.$  Using that fact, one can check that $ \overline{D_a}=\left\{D_{ax}:x\in (\F^{\times})^{2}\right\}.$ It is well known that $|(\F^{\times})^{2}|=\frac{q-1}{2}$ (See \cite{Lidl-finitefields} for example),  therefore we have $|\overline{D_a}|=\frac{q-1}{2}.$

   Let $d\in \F^{\times}-(\F^{\times})^{2},$ (We can always choose such an element), then  note that $D_{1}\nsim D_d,$ since $d=d/1\notin(\F^{\times})^2,$ therefore we have $|\overline{D_{1}}\cup \overline{D_d}|=q-1=|\F^{\times}|.$

    The last analysis implies that for any $T\in \mathcal{TM}_4(3)$ necessarily $T\sim D_1$ or $T\sim D_{d}.$ Since $D_1\nsim D_d,$ we conclude that there are exactly two isomorphism classes of matrices  $T\in \mathcal{TM}_4(3),$ such that $\mathbb{W}(T)=(3).$
\end{proof}

Note that Corollary \ref{coro1-rep-n4} and  Lemma \ref{lemma2-count-N4-q} imply that

\begin{equation}\label{eq-pre-theo-formula-n4}
    N_4(\F_q)=5+|\overline{\mathcal{TM}_4}(3,4)|.
\end{equation}

\begin{theorem}\label{theo-formula-N4}
    The number of isomorphism classes $N_4(\F_q)$ is given by:

\begin{itemize}
    \item $ N_{4}(\F_q)={(3q+23)}/{4}$ if $q\equiv 3(\textrm{mod}\quad 4).$
    \item $ N_{4}(\F_q)={(3q+25)}/{4}$ if $q\equiv 1(\textrm{mod}\quad 4).$
\end{itemize}
\end{theorem}

\begin{proof}
    By equation \ref{eq-pre-theo-formula-n4}, we only have to determine the exact value of $|\overline{\mathcal{TM}_4}(3,4)|.$ By proposition \ref{prop-representatives-n4-general}, we have $\overline{\mathcal{TM}_4}(3,4)=\{\overline{C_a}:a\in \F_q\},$  then $|\overline{\mathcal{TM}_4}(3,4)|\leq q,$ but if we use Proposition \ref{prop2-red-form-n4-W34}, we can see that  $|\overline{\mathcal{TM}_4}(3,4)|< q.$

    Let us define the following function $g:\F_q\rightarrow\F_q,$ given by $g(x)=2x+1.$ Clearly $g$ is a bijection, therefore we have:
    \begin{equation}
        |g^{-1}\left(\F_q-(\F_q^{\times})^{2}\right)|=(q+1)/2,\quad \textrm{and}\quad|g^{-1}\left((\F_q^{\times})^{2}\right)|=(q-1)/2.
    \end{equation}
     In particular note that $-1/2\in f^{-1}\left(\F_q-(\F_q^{\times})^{2}\right).$ Let us write

    \begin{equation}
        \overline{\mathcal{TM}_4}(3,4)=\left\{\overline{C_a}:g(a)\notin (\F_q^{\times})^{2}\right\}\cup\left\{\overline{C_a}:g(a)\in(\F_q^{\times})^{2}\right\},
    \end{equation}

    by Corollary \ref{coro-prop2-red-form-n4-W34}, we know that if $g(a)\notin(\F_q^{\times})^{2}$ then the condition $C_a\sim C_{a'}$ implies that $a'=a.$ Therefore
     \begin{equation}\label{eq-proof-theo-formula-N4}
        |\overline{\mathcal{TM}_4}(3,4)|=(q+1)/2+|\left\{\overline{C_a}:g(a)\in(\F_q^{\times})^{2}\right\}|
    \end{equation}

    Now consider the function $f:\F_q-\{1/2\}\rightarrow\F_q,$ given by $f(x)=-x/(2x+1).$ The number $|\left\{\overline{C_a}:g(a)\in(\F_q^{\times})^{2}\right\}|$ depends on how many fixed points of $f$ are there in the set  $g^{-1}\left((\F_q^{\times})^{2}\right).$ The fixed points of $f$ are $x=0$ and $x=-1$ and it is clear that $0\in g^{-1}\left((\F_q^{\times})^{2}\right).$

    \begin{itemize}
    \item If we assume that $-1\in g^{-1}\left((\F_q^{\times})^{2}\right),$ that is $g(-1)=-1\in(\F_q^{\times})^{2},$ then
    \begin{equation}\label{eq2-proof-theo-formula-N4}
     |\left\{\overline{C_a}:g(a)\in(\F_q^{\times})^{2}\right\}|=\frac{1}{2}\left(\frac{q-1}{2}-2\right)+2=\frac{q+3}{4}.
    \end{equation}

    It is well known that $-1\in(\F^{\times})^{2}$ if and only if $q\equiv 1(\textrm{mod}\quad 4)$ (See \cite{Lidl-finitefields} for example), therefore the number $(q+3)/4$ is an integer. Now combining equations \ref{eq-pre-theo-formula-n4}, \ref{eq-proof-theo-formula-N4} and \ref{eq2-proof-theo-formula-N4}, we obtain for this case
    \begin{equation*}
        N_4(\F_q)=5+\frac{q+1}{2}+\frac{q+3}{4}=\frac{3q+25}{4}.
    \end{equation*}
    as desired.
    \item In contrast, if we assume that $-1\notin (\F^{\times})^{2},$ then
    \begin{equation}\label{eq3-proof-theo-formula-N4}
     |\left\{\overline{C_a}:g(a)\in(\F_q^{\times})^{2}\right\}|=\frac{1}{2}\left(\frac{q-1}{2}-1\right)+1=\frac{q+1}{4}.
    \end{equation}
    Note that in this case we have that $q\equiv 3(mod\quad4),$ then the number $(q+1)/4$ in an integer. Now combining equations \ref{eq-pre-theo-formula-n4}, \ref{eq-proof-theo-formula-N4} and \ref{eq3-proof-theo-formula-N4}, we obtain for this case
    \begin{equation*}
        N_4(\F_q)=5+\frac{q+1}{2}+\frac{q+1}{4}=\frac{3q+23}{4}.
    \end{equation*}
    as desired.
    \end{itemize}
\end{proof}

\subsection{A complete set of representatives for $\mathcal{MT}_4$ over $\mathbb{C}$}\label{ssec-over-C}

In this subsection we consider as a ground field, the field of complex numbers, that is $\F=\mathbb{C}.$ If $z\in\mathbb{C},$ we denote by $|z|$ and $\mathcal{I}(z)$ its module and its imaginary part respectively. The fact that $\mathbb{C}$ is an algebraically closed field together with its geometric representation, make this case particularly interesting.

Let us denote by $\bcalU$ and $\bcalV$ the following subset of $\mathbb{C}:$

\begin{equation*}
    \bcalU=\left\{z\in \mathbb{C}:\left|z+\frac{1}{2}\right|<\frac{1}{2}\right\}\cup
    \left\{z\in \mathbb{C}:\left|z+\frac{1}{2}\right|=\frac{1}{2}\wedge \mathcal{I}(z)\geq0\right\}.
\end{equation*}

\begin{equation*}
    \bcalV=\left\{z\in \mathbb{C}:\left|z+\frac{1}{2}\right|>\frac{1}{2}\right\}\cup
    \left\{z\in \mathbb{C}:\left|z+\frac{1}{2}\right|=\frac{1}{2}\wedge \mathcal{I}(z)\leq0\right\}.
\end{equation*}

We also denote $\bcalU^{\times}=\bcalU-\{-1/2\}.$ Geometrically, in the complex plane, $\bcalU$ is the open disk of center $-1/2$ and radius $1/2,$ together with the upper half of its border, while $\bcalV$ is the exterior zone of the disk together with the lower half of its border. On the other hand, $\bcalU^{\times}$ is equal to $\bcalU$ punched at its center. Note that $\bcalU\cup\bcalV=\mathbb{C}$ and $\bcalU\cap \bcalV=\{-1,0\}.$

\begin{proposition}\label{prop1-Cz-over-C}
   \begin{enumerate}
       \item If $z\neq z'$ are two point of $\bcalU,$ then $C_{z}\nsim C_{z'}.$
       \item  If $z'\in\calV,$ then there is a unique $z\in \bcalU^{\times},$ such that $C_z\sim C_{z'}.$
   \end{enumerate}
\end{proposition}

\begin{proof}
It is enough to notice that the Möbius transformation $f:\mathbb{C}-\{-1/2\}\rightarrow \mathbb{C}$ given by $f(z)=-z/(2z+1)$ is bijective and satisfies $f(\bcalU^{\times})=\bcalV$ (See \cite{Lang-complex} for example). Also note that $f$ has two fixed points, $z=0$ and $z=-1$ and a pole in $z=-1/2,$  then if $z\in\{-1,0,-1/2\}$ we have that $C_z\sim C_{z'}$ implies that $z'=z$ (compare with Corollary \ref{coro-prop2-red-form-n4-W34}).
\end{proof}




\begin{lemma}\label{lemma-D1-over-C}
 If $\F$ is an algebraically closed field, then $|\overline{\mathcal{TM}_4}(3)|=1$
\end{lemma}

\begin{proof}
 Let $u\in \F^{\times}.$ Since $\F$ is algebraically closed, $u$ is an square on $\F,$ therefore Proposition \ref{proposition-Nn-not-finite} implies that $D_u\sim D_1,$ and therefore $\overline{\mathcal{TM}_4}(3)=\{\overline{D_1}\}$ and  $|\overline{\mathcal{TM}_4}(3)|=1.$
\end{proof}

\begin{theorem}\label{theo-full-set-for-n4-over-C}
The following set of matrices:
\begin{equation}
\left\{ 0_4,B,B',D_1 \right\}\cup \left\{C_z: z\in\bcalU\right\}.
\end{equation}
 is a complete set of representatives of (pairwise different) classes for $\mathcal{TM}_4$ over $\mathbb{C}.$
\end{theorem}

\begin{proof}
The theorem follows from Propositions \ref{prop-representatives-n4-general}, \ref{prop1-Cz-over-C} and Lemma \ref{lemma-D1-over-C}.
\end{proof}


\subsection{Isomorphisms vs ETO's for $n=4$}\label{ssec-Iso-vs-ETO-n4}
In \cite{Lobos3}, we have established the \textbf{Isomorphism vs ETO's conjecture}. That conjecture claims that for two \nSLTM, say $T$ and $S$ we have that if $T\sim S$ then $T$ and $S$ are equivalent by ETO's. As a consequence of Lemma \ref{lemma1-red-form-n4-W34} and Proposition \ref{prop2-red-form-n4-W34}, we obtain the following evidence of the veracity of the conjecture:

\begin{theorem}\label{theo-evidence-ETOs-n4}
    If $T\neq S$ are two $4-${\rm SLTM}, such that $T\sim S,$ then $T$ and $S$ are equivalent by ETO's.
\end{theorem}

\begin{proof}
Since any $4-${\rm SLTM} is equivalent by ETO's to a matrix in a simplified form (see Proposition \ref{prop-normalized-form}, we can assume that $T$ and $S$ are in simplified form. If we assume that $T\sim S,$ by Theorem \ref{theo-invariants-1REF}, we know that $\mathbb{W}(T)=\mathbb{W}(S).$ By Proposition \ref{prop-representatives-n4-general}, we only have to analyze the cases where both matrices have $\mathbb{W}=(3)$ or $\mathbb{W}=(3,4)$ as a wall.
\begin{itemize}
    \item If we assume that $\mathbb{W}(T)=\mathbb{W}(S)=(3),$ then we can assume that $T=D_a$ and $S=D_{b}$ for certain $a,b\in \F^{\times}.$ We know Proposition \ref{proposition-Nn-not-finite} , that $D_a\sim D_b$ if and only if $a/b\in(\F^{\times})^{2}.$ Assuming that, let us denote $\epsilon=\sqrt{a/b}$ and define progressively the following matrices
    \begin{equation}
    \begin{array}{ccc}
        T_2:=\calQ_{(2,1)}(D_a,\epsilon), &  T_3:=\calQ_{(3,2)}\left(T_2,{(\epsilon-1)/(2\epsilon)}\right), &  T_4:=\calP_{2}\left(T_3,\epsilon\right)\\
         \quad& \quad &\quad\\
         T_5:=\calQ_{(2,1)}(T_4,-1)&T_6:=\calP_{1}(T_5,1/\epsilon)&T_7:=\calP_{4}(T_6,\epsilon)
    \end{array}
    \end{equation}
    By construction we have $D_a$ is equivalent by ETO's to $T_7.$ If we do the calculus, we could see that $T_7=D_b,$ we left to the readers to check this.

    \item If we assume that $\mathbb{W}(T)=\mathbb{W}(S)=(3,4),$ then we can assume that $T=C_a$ and $S=C_{a'}$ for certain $a,a'\in \F.$ Since $T\sim S,$ Proposition \ref{prop2-red-form-n4-W34}, implies that $2a+1\in(\F^{\times})^{2}$ and $a'=-\frac{a}{2a+1}.$  Under this conditions, the proof of Proposition \ref{prop2-red-form-n4-W34}, implies that $T\approx S.$
\end{itemize}

\end{proof}


\begin{appendices}


\section{Notations and terminology}\label{ssec-Notations}
\begin{itemize}
  \item $\F$ is a fixed field of characteristic different from $2.$ We denote $\F^{\times}=\F-\{0\},$ $\F^{2}=\{a^{2}:a\in \F\}$ and $(\F^{\times})^{2}=\F^{2}-\{0\}.$
  \item A \nSLTM,  is a strictly lower $n\times n-$matrix,  with entries in $\F.$ A \SLTM, is a strictly lower matrix (of any size), with entries in $\F.$ We denote by \nSLTMn, the set of all \nSLTM,  and by \SLTMn,  the set of all \SLTM.
  \item Let $a,b\in\mathbb{Z}$, with $a\leq b.$ we define the interval $[a,b]:=\{m\in\mathbb{Z}:a\leq m \leq b\}.$ Analogously we define $[a,b[,]a,b[,$ etc. We also denote $[a]=[a,a]=\{a\}.$
  \item If $T$ is a \nSLTM,  $I,J\subset [1,n]$ two intervals, then we denote by $T^{I}_{J}$ the submatrix of $T$ comprised by the entries $t_{ij}$ where $i\in I$ and $j\in J.$ We also denote $T^{I}:=T^{I}_{[1,n]}$ and $T_{J}:=T^{[1,n]}_{J}.$ In  particular $T^{[r]}$ denotes the $r-$th row of $T$ and $T_{[c]}$ denotes the $c-$th column of $T$
  \item If $T$ is a \nSLTM, then for each $r\in[1,n]$ we define:
  \begin{equation*}
    c_r=\left\{\begin{array}{cc}
                 \max\{j:t_{rj}\neq 0\}, & \textrm{if}\quad T^{[r]}\neq 0. \\
                 0, & \textrm{otherwise}
               \end{array}\right.
  \end{equation*}
  When $c_r>0$ we say that the entry $t_{r,c_r}$ is the \emph{leader} of the $r-$th row of $T.$ The pair $(r,c_r)$ is the \emph{leader position} of the $r-$th row of $T.$
  \item Given a \nSLTM, $T=[t_{ij}]$ and $\alpha\in\F.$ For each $1\leq i<j<k\leq n,$ we denote $\Delta_{i,j,k}^{(\alpha)}(T)=\alpha t_{ki}+t_{kj}t_{ji}.$ If there is no possible confusion we only write $\Delta_{i,j,k}^{(\alpha)}.$
  \item Given any matrix $M=[m_{ij}],$ we denote by $\mu(M)=|\{(i,j):m_{ij}\neq 0\}|.$ We call $\mu(M)$ the \emph{measure} of the matrix $M.$
\end{itemize}
\section{Previous definitions}\label{ssec-prev-def}

This section is a compilation of definitions coming from \cite{Lobos2}, \cite{Lobos3} and \cite{Lobos4}.

\begin{definition}\label{def-A(T)}
Fix a \emph{ground field} $\F.$ For each strictly lower triangular matrix $T=[t_{ij}]$ of size $n\times n,$ with $n\geq 2$ and entries on $\F,$ we define $\calA(T)$ as $\F-$algebra given by generators $X_1,\dots,X_n$ and relations given by
 \begin{equation}\label{EQ-def-A(T)}
    \left\{\begin{array}{c}
      X_1^{2}=0 \\
       \quad\\
      X_i^{2}=\sum_{j<i}{t_{ij}}X_jX_i,\quad (2\leq i \leq n)\\
      \quad\\
      X_iX_j-X_jX_i=0,\quad (1\leq i,j \leq n)
    \end{array}\right.
  \end{equation}
 The assignment $\deg(X_i)=2,$ for $1\leq i\leq n,$ converts $\calA(T)$ in a commutative graded algebra.
 \end{definition}

\begin{definition}\label{def-category-bcalNT}
  We define the category $\bcalNT$  as the one whose objects are the many algebras  $\calA(T),$ where $T$ is a \SLTM,  with entries in the ground field $\F,$ and whose morphisms are all the \emph{preserving degree homomorphisms} between them.
\end{definition}

We write $T\sim S,$ whenever the algebras $\calA(T)$ and $\calA(S)$ are isomorphic.

\begin{definition}\label{def-ETOs}
Let $T=[t_{rk}]$ be a \nSLTM:
\begin{enumerate}
  \item  Fix an $1\leq r_1\leq n$ and an scalar $\alpha\neq 0.$ We define the matrix $S=\calP_{r_1}(T,\alpha),$ as the \nSLTM, whose entries are given by equation \ref{eq-def-P-operation}
  \begin{equation}\label{eq-def-P-operation}
    s_{rk}=\left\{\begin{array}{cc}
                    \alpha^{-1} t_{rk} & \textrm{if}\quad r=r_1 \\
                    \alpha t_{rk} & \textrm{if}\quad k=r_1 \\
                    t_{rk} & \textrm{otherwise}
                  \end{array}\right.
  \end{equation}
  \item Fix $1\leq r_1<r_2\leq n.$ If the matrix $T$ satisfies the following conditions:
  \begin{enumerate}
    \item The entry $t_{r_2j}=0$ for each $r_1\leq j$
    \item The entry $t_{rr_1}=0$ for each $r_1\leq r\leq r_2$
  \end{enumerate}
  then we define the matrix $S=\calF_{(r_1,r_2)}(T)$ as the \nSLTM, whose entries are given by the following equation:
  \begin{equation}\label{eq-def-F-operation}
    s_{rk}=\left\{\begin{array}{cc}
                    t_{r_2k} & \textrm{if}\quad r=r_1 \\
                    t_{r_1k} & \textrm{if}\quad r=r_2  \\
                    t_{rr_2} & \textrm{if}\quad k=r_1 \\
                    t_{rr_1} & \textrm{if}\quad k=r_2 \\
                    t_{rk} & \textrm{otherwise}
                  \end{array}\right.
  \end{equation}
  \item Fix a leader position $(r_0,k_0)$ of $T$
  \begin{enumerate}
    \item If $k_0=1,$ then for any $\beta\in \F^{\times},$ we define the matrix $S=\calQ_{(r_0,1)}(T,\beta)$ as the \nSLTM, whose entries $s_{rk}$ are given by equation \ref{eq-def-Q-operation}.
    \item If $k_0>1$ and there is a $\beta\in \F^{\times},$ such that $\Delta^{(1)}_{i,k_0,r_0}(T)=\beta t_{k_0,i},$ for each $1\leq i<k_0,$ then we define the matrix $S=\calQ_{(r_0,k_0)}(T,\beta)$ as the \nSLTM, whose entries $s_{rk}$ are given by equation \ref{eq-def-Q-operation}.
  \end{enumerate}
  \begin{equation}\label{eq-def-Q-operation}
    s_{rk}=\left\{\begin{array}{cc}
                    t_{r_0,k_0}-2\beta & \textrm{if}\quad (r,k)=(r_0,k_0) \\
                    \quad & \quad\\
                    t_{r,k_{0}}+\beta t_{r,r_{0}} & \textrm{if}\quad r>r_0\quad \textrm{and}\quad k=k_0  \\
                    \quad & \quad\\
                    t_{rk} & \textrm{otherwise}
                  \end{array}\right.
  \end{equation}
\end{enumerate}
We call \emph{Elementary triangular operations} (ETO for short) any of the following matrix assignments (whenever they are well defined):

\begin{equation}\label{eq-ETO-assignments}
  \begin{array}{cc}
    \calP_r(\quad,\alpha):T\mapsto \calP_r(T,\alpha)&(\calP-\textrm{operation} ) \\
    \quad&\quad  \\
     \calF_{(r_1,r_2)}(\quad):T\mapsto \calF_{(r_1,r_2)}(T)& (\calF-\textrm{operation}) \\
    \quad &\quad \\
    \calQ_{(r_0,k_0)}(\quad,\beta):T\mapsto \calQ_{(r_0,k_0)}(T,\beta)& (\calQ-\textrm{operation})
  \end{array}
\end{equation}
\end{definition}
 Let $T,S$ be two \nSLTM, if there exists of a (finite) sequence of ETO's that transform $T$ into $S,$ then we say that $T$ and $S$ are \emph{equivalent  by ETO's}, and write $T\approx S.$




\begin{definition}\label{def-first-red-form}
  Let $S\neq 0_n$ be a \nSLTM, we say that $S$ is in the \emph{First Reduced Echelon Form} ($1-$REF), if it satisfies the following conditions:
  \begin{enumerate}
    \item The zero rows of $S$ are above of the nonzero rows.
    \item For each leader position $(r,c)$ we have that $c>1$ and there is a $i_0\in [1,c[$ such that $\Delta^{(2)}_{i_0,c,r}\neq 0.$
    \item If $(r,c_r)$ and $(r+1,c_{r+1})$ are the leader positions of the $r-$th and $(r+1)-$th rows of $S,$ then $c_r\leq c_{r+1}.$
  \end{enumerate}

  We also say that the zero matrix, $0_n$ is in the first reduced echelon form.
\end{definition}

\begin{definition}\label{def-Blocks-1REF}
  Let $S$ be a \nSLTM, such that $S\nsim 0_n.$ Assuming that $S$ is in the $1-$REF:
  \begin{enumerate}
   \item We define recursively numbers $r_j$ and a vector $\mathbb{W}(S)$ associated to $S$ as follows:
      \begin{enumerate}
      \item Let $r_1:=\min\{r\in[1,n]:S^{[r]}\neq 0\}$.
      \item Assuming the number $r_j$ is already defined, define the set $A_j=\{r\in]r_j,n]:c_r>r_j\}.$
        \begin{itemize}
          \item If $A_j\neq \emptyset$ then we define $r_{j+1}=\min(A_j).$
          \item If $A_j=\emptyset,$ then we define the  vector $\mathbb{W}(S)=(r_1,\dots,r_j)$.
        \end{itemize}
        The vector $\mathbb{W}(S)$ will be called the \emph{wall} of $S,$
      \end{enumerate}
   \item If $\mathbb{W}(S)=(r_1),$ then we define the \emph{brick} of $S$ by:
     \begin{equation*}
       B_1(S):=S_{[1,r_1[}^{[r_1,n]}
     \end{equation*}
   \item If $\mathbb{W}(S)=(r_1,\dots,r_h),$ with $h>1,$ then for each $j\in[1,h]$ we define the $j-$th \emph{brick} of $S$ as the submatrix:
       \begin{equation*}
         B_{j}(S)=\left\{\begin{array}{cc}
                        S_{[1,r_1[}^{[r_1,r_{2}[} & \textrm{if}\quad j=1 \\
                        \quad & \quad \\
                        S_{[r_{j-1},r_j[}^{[r_j,r_{j+1}[} & \textrm{if}\quad 1<j<h\\
                        \quad & \quad \\
                        S_{[r_{h-1},r_h[}^{[r_h,n]}&\textrm{if}\quad j=h
                      \end{array}\right.
       \end{equation*}
  \end{enumerate}
  For the zero matrix we also define $\mathbb{W}(0_n)=(0).$ 
\end{definition}

For a given $\mathbb{W}(S)=(r_1,\dots,r_{h}),$ we will denote $r_0:=1$ and $r_{h+1}=n+1,$ then  $[r_0,r_1[=[1,r_1[$ and $[r_h,r_{h+1}[=[r_h,n+1[=[r_h,n].$

\begin{definition}\label{def-general-wall}
Let $T$ be any \SLTM. Suppose that $S$ is a \SLTM, that is in $1-$REF and $T\sim S.$ We define the \emph{wall} of $T$ by $\mathbb{W}(T)=\mathbb{W}(S).$
\end{definition}

From \cite{Lobos4}, we know that for each strictly lower triangular matrix $T,$ the wall $\mathbb{W}(T),$ is well defined.

\begin{definition}\label{def-second-red-form}
 Let $T\nsim 0_n$ be a \nSLTM, and let $\mathbb{W}(T)=(r_1,\dots,r_h)$ its wall. We say that $T$ is in the \emph{Second Reduced Echelon Form} ($2-$REF) if it satisfies the following conditions:
 \begin{enumerate}
   \item $T$ contains exactly $r_1-1$ zero rows and they are above of the nonzero rows.
   \item If $(r,c_r)$ is the leader position of the $r-$th row of $T,$ then $c_r>1$ and there is a $i_0\in [1,c_r[$ such that $\Delta^{(2)}_{i_0,c,r}\neq 0.$ Moreover for each $j\in [1,h]$ and $r\in[r_j,r_{j+1}[$ we have that, $c_r\in [r_{j-1},r_j[.$
   \item For each $j\in [1,h]$ and $r\in[r_j,r_{j+1}-1[$ we have that
   \begin{equation*}
     \mu\left(T_{[r_j,r_{j+1}[}^{[r]}\right)\leq \mu\left(T_{[r_j,r_{j+1}[}^{[r+1]}\right).
   \end{equation*}
   We also say that the zero matrix is in the $2-$REF.
 \end{enumerate}
\end{definition}

\begin{definition}\label{def-Blocks-2REF}
  Let $T\nsim 0_n$ be a \SLTM, and assume that is is in $2-$REF, let $\mathbb{W}(T)=(r_1,\dots,r_h).$ For each $j\in[1,h]$ we define recursively numbers $r_{jk}$ and $\mu_{jk}$ and a matrix $M_j(T),$ as follows:
  \begin{enumerate}
    \item Let $r_{j1}:=r_j,$ and $\mu_{j1}=\mu\left(T_{[r_{j-1},r_j[}^{[r_{j1}]}\right).$
    \item Assuming that the numbers $r_{jk}$ and $\mu_{jk}$ are already defined. Define the set
         \begin{equation*}
           A_{jk}=\left\{r\in]r_{jk},r_{j+1}[\quad : \quad \mu\left(T_{[r_{j-1},r_j[}^{[r]}\right)>\mu_{jk}\right\}.
         \end{equation*}
         \begin{itemize}
           \item If $A_{jk}\neq \emptyset,$ then define $r_{j,k+1}:=\min(A_{jk})$ and $\mu_{j,k+1}:=\mu\left(T_{[r_{j-1},r_j[}^{[r_{j,k+1}]}\right)$
           \item If $A_{jk}=\emptyset,$ then define the matrix $M_j(T)$ as follows:
               \begin{equation*}
                 M_j(T):=\left[\begin{matrix}
                                 r_{j1} & r_{j2} & \cdots & r_{jk} \\
                                 \mu_{j1} & \mu_{j2} & \cdots & \mu_{jk}
                               \end{matrix}\right].
               \end{equation*}
         \end{itemize}
          We call the sequence $\mathbb{M}(T):=(M_1(T),\dots,M_h(T))$ the \emph{measure sequence} of $T.$
  \end{enumerate}
\end{definition}

\begin{definition}
Let $T$ be a \SLTM. We define the \emph{measure sequence} of $T$ by $\mathbb{M}(T)=\mathbb{M}(V),$ where $V$ is any matrix in $2-$REF, such that $T\sim V.$
\end{definition}

From \cite{Lobos4}, we know that for each strictly lower triangular matrix $T\nsim 0_n,$ the measure sequence $\mathbb{M}(T),$ is well defined.

\section{Previous results}\label{ssec-prev-results}

This section is a compilation of previous results coming from \cite{Lobos2}, \cite{Lobos3} and \cite{Lobos4}.

\begin{theorem}\label{theo-monomial-basis}
  If $T=[t_{ij}]$ is a \nSLTM, then the dimension of the algebra $\calA(T)$ is $2^n$ and the set $\{X_1^{a_1}\cdots X_n^{a_n}:a_i\in\{0,1\}\}$ is a basis for $\calA(T).$
\end{theorem}
\begin{proof}
  See \cite{Lobos2}.
\end{proof}

\begin{theorem}\label{theo-Key-condition}
Let $T,S$ be two \nSLTM.  Given a $n\times n$ invertible matrix, $\Gamma=[\gamma_{ij}],$ then the assignment
  \begin{equation*}
    X_j\mapsto \sum_{i=1}^{n}{\gamma_{ij}}{Y_i}
  \end{equation*}
  defines an isomorphism $\bgamma:\calA(T)\rightarrow\calA(S),$ if and only if the following system of equations holds:
  \begin{equation}\label{Eq1-theo-Key-condition}
    2\gamma_{ir}\gamma_{kr}+\gamma_{kr}^{2}s_{ki}=
    \sum_{j<r}{t_{rj}}(\gamma_{kj}\gamma_{kr}s_{ki}
    +\gamma_{kj}\gamma_{ir}+\gamma_{ij}\gamma_{kr})
    ,\quad (1\leq r\leq n;1\leq i<k\leq n).
  \end{equation}
\end{theorem}
\begin{proof}
  See \cite{Lobos3}.
\end{proof}

\begin{theorem}\label{theo-sim-vs-approx}
  Let $T,S$ be two \nSLTM.  If there exists a (finite) sequence of ETO's that transform $T$ into $S,$ then $T\sim S.$
\end{theorem}
\begin{proof}
  See \cite{Lobos3}.
\end{proof}

 \begin{proposition}\label{Prop-Agorithm-1REF}
  Let $T$ be a \nSLTM, such that $T\nsim 0_n.$ Then $T$ is equivalent by ETO's to a  matrix $R$ that is in the first reduced echelon form.
\end{proposition}

\begin{proof}
  See \cite{Lobos4}
\end{proof}

\begin{proposition}\label{Prop-Agorithm-2REF}
  Let $T$ be a \nSLTM, such that $T\nsim 0_n.$ Then $T$ is equivalent by ETO's to a  matrix $R$ that is in the second reduced echelon form.
\end{proposition}

\begin{proof}
  See \cite{Lobos4}
\end{proof}

 \begin{theorem}\label{theo-invariants-1REF} $(\boldsymbol{\mathbb{W}-}\textrm{\textbf{Test}})$
  Let $S,T$ be two \nSLTM, If $S\sim T,$ then $\mathbb{W}(S)=\mathbb{W}(T).$
\end{theorem}

\begin{proof}
  See \cite{Lobos4}
\end{proof}

\begin{proposition}\label{prop-Gamma-triang-2}
  Let $T,S$ be two nonzero \nSLTM, such that:
   \begin{itemize}
     \item $T$ and $S$ are in the $2-$REF,
     \item $\mathbb{W}(T)=\mathbb{W}(S)=(r_1,\dots,r_h).$
     \item $\mathbb{M}(T)=\mathbb{M}(S)=(M_1,\dots,M_h).$
   \end{itemize}

     If $\bgamma:\calA(T)\rightarrow\calA(S)$ is an isomorphism, then the associated matrix $\Gamma=[\gamma_{ij}]$ satisfies the following relations:
  \begin{enumerate}
    \item If $r\in [r_j, r_{j+1}[$ then $\gamma_{kr}=0$ for all $k\geq r_{j+1}.$
    \item There is a bijection $\rrn:[1,n]\rightarrow [1,n]$ such that:
    \begin{enumerate}
      \item For each $j\in[1,h]$ we have $\rrn([r_j,r_{j+1}[)= [r_j,r_{j+1}[$
      \item For each $i,k\in [r_j,r_{j+1}[$ we have that
    \begin{equation*}
       \gamma_{ik}\neq 0 \quad \textrm{if and only if}\quad k=\rrn(i).
    \end{equation*}
    \item For each $j\in[1,h],$ if
    \begin{equation*}
      M_j=\left[\begin{matrix}
                     r_{j1} & \cdots & r_{j,l(j)} \\
                     \mu_{j1} & \cdots  & \mu_{j,l(j)}
                   \end{matrix}\right]
    \end{equation*}
    for some $l=l(j)>1.$  Then for each $i\in[1,l(j)[$ we have that $\rrn([r_{j,i},r_{j,i+1}[)= [r_{j,i},r_{j,i+1}[.$
    \end{enumerate}
  \end{enumerate}
\end{proposition}

\begin{proof}
See \cite{Lobos4}.
\end{proof}

\begin{theorem}\label{theo-invariants-2REF} $(\boldsymbol{\mathbb{M}-}\textrm{\textbf{Test}})$
 Let $T,S$ be two \nSLTM. If $T\sim S,$ then $\mathbb{M}(T)=\mathbb{M}(S).$
\end{theorem}
\begin{proof}
 See \cite{Lobos4}.
\end{proof}






\end{appendices}

\sc
\begin{itemize}
\item diego.lobosm@uv.cl, {Universidad de Valpara\'{\i}so, Chile.}
\end{itemize}

\end{document}